\documentclass[12pt]{article}
\usepackage{amsmath}
\usepackage{amssymb}
\usepackage{amsfonts}
\usepackage{eucal}
\usepackage[usenames]{color}
\usepackage{graphicx}
\numberwithin{equation}{section}
\newtheorem{Theorem}{Theorem}[section]
\newtheorem{Proposition}[Theorem]{Proposition}

\newtheorem{Remark}[Theorem]{Remark}
\begin{document}
\title{On the joint distribution of first-passage time and first-passage area of  drifted Brownian motion}
\author{Mario Abundo$^*$, Danilo Del Vescovo\thanks{Dipartimento di Matematica, Universit\`a  ``Tor Vergata'', via della Ricerca Scientifica, I-00133 Rome, Italy.
Corresponding author E-mail: \tt{abundo@mat.uniroma2.it}}
}
\date{}
\maketitle

\begin{abstract}
\noindent
For drifted Brownian motion $X(t)= x - \mu t + B_t \ (\mu >0)$ starting from $x>0,$ we study the joint distribution
of the first-passage time below zero, $\tau(x),$ and
the first-passage area, $A(x),$  swept out by $X$ till the time $\tau(x).$ In particular, we establish differential equations
with boundary conditions for the joint moments $E[\tau(x)^m A(x)^n],$ and we present an  algorithm  to find recursively them,
for any $m$ and $n.$ Finally, the expected value of the time average of $X$ till  the time $\tau(x)$ is obtained.
\end{abstract}


\section{Introduction}
This a continuation of the article \cite{abundo:mcap13}, where
we have studied the probability distribution
of the first-passage area (FPA)  $A(x),$ swept out by a one-dimensional jump-diffusion process $X(t),$ starting from $x>0,$
till its first-passage time (FPT) $\tau(x)$
below zero. In particular, in \cite{abundo:mcap13} it was shown that the Laplace transform
and the moments of $A(x)$ are solutions to certain partial
differential-difference equations with outer conditions. \par
In the present paper, we aim to study the joint distribution of $\tau(x)$ and $A(x),$ in the special case when
$X(t)$ is Brownian motion with negative drift $- \mu,$ that is, $X(t)= x - \mu t + B_t, $ where $B_t$ denotes standard Brownian
motion.
We state differential equations for the Laplace transform
of the two-dimensional random variable $(\tau(x), A(x)),$ and for
the joint moments  $E[\tau(x)^m A(x)^n]$ of the FPT and FPA; moreover, we present an  algorithm  to find them recursively,
for any $m$ and $n.$
In particular, closed formulae for $E[\tau(x) A(x)]$ and   the correlation between $\tau(x)$ and $A(x),$ are
explicitly obtained; they match existing results, obtained in \cite{kearney:jph14} in terms of the Airy function $Ai(z)$ by
using its asymptotic expansion as $z \rightarrow \infty. $
Furthermore, we find the expected value of the time average of $X(t)$ till its FPT below zero. \par\noindent
The difference compared to similar articles (see e.g. \cite{kearney:jph14}) is that the results of the present paper are obtained
without the use of special functions. \par
Over the years, several authors
have dealt with first-passage area  for Brownian
motion with negative drift
(see \cite{janson:pro07}, \cite{kearney:jph14}, \cite{kearney:jph05}, \cite{kearney:jph07},
\cite{knight:jam00}, \cite{perman:aap96} and references in \cite{abundo:mcap13});  the FPA
has interesting applications, for instance, in Queueing Theory, if one
identifies $X(t)$ with the length of a queue at time $t,$ and
$\tau(x)$ with the busy period, that is the time until the queue
is first empty; then, $A(x)$ represents the cumulative waiting
time experienced by all the ``customers'' during a busy period.
The FPA arises also in solar physics studies, non-oriented animal movement patterns, and
DNA breathing dynamics (see references in \cite{kearney:jph14}). For other examples from Economics and Biology, see \cite{abundo:mcap13}.
\par
The paper is organized as follows. Section 2 contains the main results: in subsection 2.1 we deal
with the Laplace transform
of $(\tau(x), A(x))$ and the joint moments of $A(x)$ and $\tau (x),$ precisely, in 2.1.1
we find explicitly $E[\tau(x) A(x)],$ and in 2.1.2 we establish ODEs for
the joint moments
 $E[\tau(x)^m A(x)^n]$ of order $n+m, \ n, m \ge 0;$  moreover we present
an algorithm to find recursively $E[\tau(x)^m A(x)^n],$ for any $m$ and $n;$
subsection 2.2 is devoted to  find the expected value of the time average of $X(t)$ till its FPT below zero.
Section 3 contains conclusions and final remarks.

\section{Notations, formulation of the problem and main results}
For $x>0$ and $\mu >0,$  we consider Brownian motion with negative drift $- \mu, $  that is $X(t)= x - \mu t + B_t, $ where $B_t$ is standard Brownian
motion.
We denote by
\begin{equation} \label{firstpassagetime}
\tau(x) = \tau ^ \mu (x)= \inf \{ t >0: X(t) \le 0 | X(0)=x \}
\end{equation}
the FPT below zero of $X(t);$  since the drift $- \mu $ is negative, $\tau(x)$ is finite with
probability one.
Let $U$ be a functional of the process $X,$ and for $\lambda >0$ denote by
\begin{equation}
M_ {U, \lambda} (x) =  E \left [ e^ { - \lambda
\int _0^ {\tau(x)} U(X(s)) ds} \right ]
\end{equation}
the Laplace transform of the integral $\int _0 ^ {\tau(x)}
U(X(s))ds.$ We recall from \cite{abundo:mcap13} the following:
\begin{Theorem}
For $\lambda >0, \ M_{ U, \lambda} (x) =  E \left [ e^ { - \lambda
\int _0^ {\tau(x)} U(X(s)) ds} \right ]$ is the solution of the problem with boundary conditions ($M'$ and $M''$ denote first and second derivative with
respect to $x)$ :
\begin{equation} \label{brownianlaplace}
\begin{cases}
\frac 1 2 M'' _{U, \lambda } (x) - \mu  M'_ {U, \lambda } (x) =
\lambda U(x) M_ {U, \lambda } (x ) \\
M_ {U, \lambda }(0) =1 \\
\lim _ { x \rightarrow + \infty} M_ {U, \lambda } (x) =0 \ .
\end{cases}
\end{equation}
Moreover,   the n-th order moments $T_n(x)$ of $ \int _0 ^ {\tau(x)}
U(X(s))ds,$ if they exist finite,
satisfy the recursive ODEs:
\begin{equation} \label{eqmomentsimple}
\frac 1 2  T''_n(x) - \mu T'_n(x) = -n U(x) T_{n-1} (x), \ {\rm for} \  x >0
\end{equation}
with the condition $T_n(0)= 0,$ plus an  appropriate  additional  condition.
\end{Theorem}
\hfill $\Box$
\bigskip

\noindent
{\bf (i)} The moment generating function  of $\tau ^ \mu(x)$
\par \noindent By solving \eqref{brownianlaplace}  with $U(x)=1,$ one
explicitly obtains:
\begin{equation} \label{Laplacetau}
M_ {U, \lambda } (x)= E(\exp( -\lambda \tau ^ \mu (x))= \exp [-x( \sqrt { \mu ^2 + 2 \lambda } - \mu )] .
\end{equation}
This Laplace transform can be explicitly inverted (see  \cite{kartay:sto75}), so obtaining the well-known expression of the density of $\tau ^\mu(x):$
\begin{equation} \label{brownianfirstpassagedensity}
f_ {\tau ^\mu(x)} (t) = \frac {x } {\sqrt { 2 \pi } t^ { 3/2} } e^ {-(x - \mu t)^2/ 2t } ,
\end{equation}
that is the inverse Gaussian density.
For $\mu >0$ the moments $T_n (x)= E(\tau^ \mu(x))^n$ of any order $n,$ are finite and they can be easily obtained by calculating
$ (-1)^n  [\partial ^n  M _ {U, \lambda } (x) / \partial \lambda ^n]_ { \lambda =0}.$
We obtain, for instance, $E(\tau ^\mu(x))= \frac x \mu $ and $E( (\tau^\mu (x))^2)= \frac x { \mu ^3} + \frac {x^2} { \mu ^2} .$
Note that $E(\tau ^\mu(x)) \rightarrow + \infty,$ as $ \mu \rightarrow 0.$
As easily seen, for any $x>0,$ it results $M_ {U, \lambda } (x) \rightarrow 1,$ as $\mu \rightarrow + \infty,$ or, equivalently,
$\tau ^\mu(x)$ converges to $0$ in distribution, and so $E((\tau ^\mu(x))^n) \rightarrow 0,$ as $\mu \rightarrow + \infty,$ for any $n$ and $x>0.$

\noindent {\bf (ii)} The moments of $A^\mu(x)= \int _0 ^ { \tau ^ \mu (x)} (x- \mu t + B_t) dt .$ \par
\noindent For $U(x)=x$ the equation of problem \eqref{brownianlaplace} becomes
\begin{equation} \label{eqlaplacearea}
\frac 1 2 M''_ {U, \lambda } - \mu M' _ {U, \lambda } = \lambda x M_ {U, \lambda } ,
\end{equation}
where now $M_{ U, \lambda} (x) = E( e ^ { - \lambda A^ \mu (x) } ) .$
It is  the Schrodinger equation for a quantum particle moving in a
uniform field (see e.g \cite{kearney:jph05},
\cite{kearney:jph07}). Its explicit solution
cannot be found in terms of elementary functions, however it can be
written in terms of the Airy function (see
\cite{grand:tab80}, \cite{kearney:jph05}) though, for $\mu \neq 0,$ it is impossible to
invert the Laplace transform $M_ {U, \lambda }$  to obtain the
probability density of $A ^ \mu (x).$ \par\noindent In the special case
$\mu =0,$ it can be shown (\cite{kearney:jph05}) that the solution
of \eqref{eqlaplacearea} is:
\begin{equation}
M_{ U, \lambda} (x) = 3 ^ {2/3} \Gamma \left ( \frac 2 3 \right ) {\rm
Ai} (2 ^ {1/3} \lambda ^ {1/3} x ),
\end{equation}
where ${\rm Ai }(z)$ denotes the standard Airy function, which obeys the ODE ${\rm Ai} '' (z) = z {\rm Ai} (z)$ (see \cite{grand:tab80}); then, by  inverting
this Laplace transform one finds that the first-passage area
density is (\cite{kearney:jph05}):
\begin{equation}
f_ {A^ 0 (x)} (a) = \frac {2^ {1/3} } {3 ^ {2/3} \Gamma ( \frac 1 3 )
} \ \frac {x } {a ^ {4/3} } \ e^ { -2 x^3 / 9a } \ .
\end{equation}
Thus,  the distribution of $A^ 0(x)$ has an algebraic
tail of order $ \frac 4 3 $ and so the moments of all orders are
infinite. \par
For $\mu >0,$
closed form expression for the first two moments of $A ^ \mu(x)$ have been obtained in \cite{abundo:mcap13},
by solving \eqref{eqmomentsimple} with $U(x) =x ,$ with appropriate conditions.
The mean first-passage area is
\begin{equation} \label{expressionmeanarea}
E(A ^ \mu (x))= \frac {x^2 } {2 \mu } + \frac {x } {2 \mu ^2 } \ ,
\end{equation}
while the second moment of $A ^ \mu (x)$ is
\begin{equation} \label{expressionsecondarea}
E[(A^\mu (x))^2]= \frac {x^4 } {4 \mu ^2 } + \frac {5 x^3 } {6 \mu ^3 } + \frac {5x^2 } {4 \mu ^4 } + \frac {5 x } {4 \mu ^5 } \ .
\end{equation}
Finally, by \eqref{expressionmeanarea} and \eqref{expressionsecondarea}, the variance of $A^ \mu (x)$ is obtained:
\begin{equation} \label{variancearea}
Var(A ^\mu (x)) = E (A^ \mu (x) ^2 ) - E^2 (A ^ \mu (x)) = \frac {x^3 } {3 \mu ^3 } + \frac {x^2 } { \mu ^4 } + \frac {5x } {4 \mu ^5 } \ .
\end{equation}
Since a closed form expression of the density of the first-passage area $A^ \mu (x)$ cannot be found for $\mu >0 ,$
it must be studied numerically;
indeed, it was estimated 
in \cite{abundo:mcap13} by simulating a large number of trajectories of Brownian motion
with drift $- \mu,$ starting from the initial state $x>0.$ For completeness, we report Figure 1 from \cite{abundo:mcap13}, which shows the shape of the
estimated FPA density for $x=1,$  as a function of $\mu >0.$

\begin{figure}
\centering
\includegraphics[height=.4\textheight]{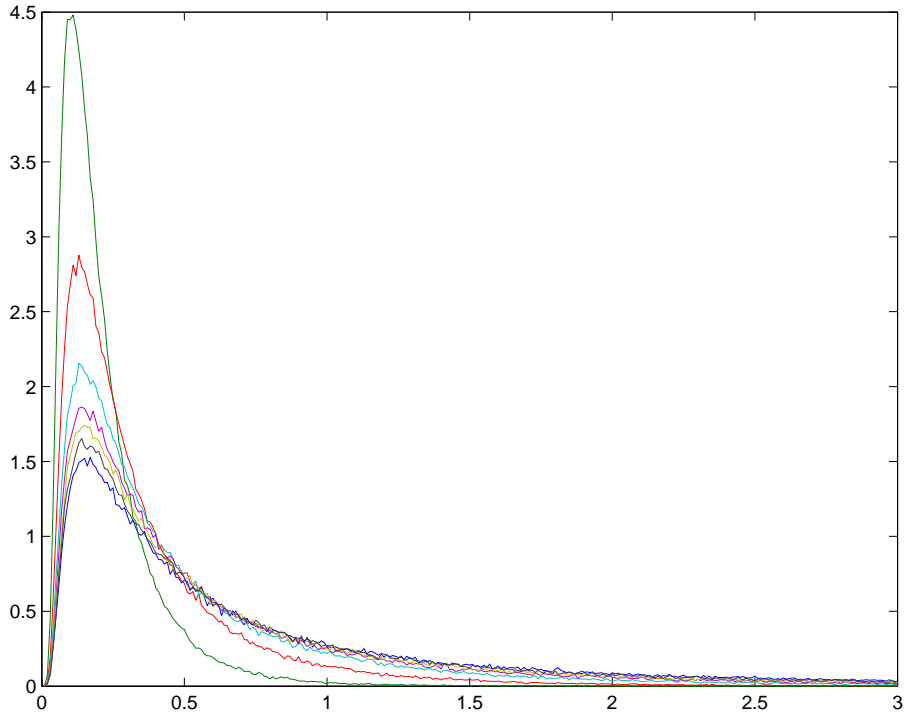}
\caption{Estimated density of the first-passage area
$A^ \mu (x)$ for $x=1$ and several values of $\mu.$ From  top to bottom,
with respect to the peak of the curve: $\mu = 3; \ \mu = 2;  \ \mu
= 1.5; \ \mu = 1.3; \  \mu = 1.2 ; \  \mu = 1.1 ; \ \mu = 1.$
 }
\end{figure}

\subsection{Joint moments of $A(x)$ and $\tau (x).$}
Now, we go to study the joint distribution of $A(x)$ and $\tau (x);$ our aim is to find an explicit form for
the joint moments of $A(x)$ and $\tau (x).$
For $\lambda = (\lambda_1, \lambda _2 ) \in \mathbb{R}_+ ^ 2 ,$ the joint Laplace transform of $(A(x), \tau (x))$ is
$$M_ \lambda (x) = M_ {(\lambda _1, \lambda _2)} (x) =E \left ( e^ { - \lambda _1 \tau (x) - \lambda _2 A(x) } \right ) $$
\begin{equation} \label{joint Laplace}
= E \left ( e^ { - \lambda _1 \tau (x) - \lambda _2 \int _0 ^ { \tau (x) }  X(t) dt } \right ) =
E \left ( e^ { - \int _0 ^ {\tau(x)} (\lambda _1 + \lambda _2 X(t) ) dt }\right ).
\end{equation}
By using \eqref{brownianlaplace} with $\lambda =1$ and
$U(x)= \lambda _1 + \lambda _2 x,$  we obtain that the function
$M_ {(\lambda _1, \lambda _2 )} (x)$ solves the problem:
\begin{equation} \label{jointLaplaceeq}
\begin{cases}
\frac 1 2 \frac {d^2} {dx^2} M _ \lambda (x) - \mu \frac d {dx} M _ \lambda (x) = ( \lambda _1 + \lambda _2 x) M _ \lambda (x) , \\
M _ \lambda (0) =1, \\ \lim _ { x \rightarrow + \infty } M _ \lambda (x)=0 .
\end{cases}
\end{equation}
Explicitly one has (see \cite{kearney:jph14}):
\begin{equation} \label{eqforjointLaplace}
M _ \lambda (x) = e^{\mu x }
\left [ Ai \left ( \frac {\mu ^2 + 2 \lambda _1 + 2 \lambda _2 x } {2^{2/3} \lambda _2 ^{2/3} }  \right ) \right ] /
\left [ Ai \left ( \frac {\mu ^2 + 2 \lambda _1  } {2^{2/3} \lambda _2 ^{2/3} }\right ) \right ],
\end{equation}
where ${\rm Ai} (z)$ is the standard Airy function. The joint moments $E [ \tau ^m (x) A^n (x)],$ if they exist  finite,
can be obtained by taking the partial derivative of order $m+n$ at $\lambda = 0,$  i.e.
$$ E \left [ \tau ^m (x) A^n (x) \right ] =
(-1)^{m+n} \frac {\partial ^m } { \partial \lambda _1 ^m} \frac {\partial ^n } { \partial \lambda _2 ^n} M_ \lambda (x) \Big |
_{\lambda _1 = \lambda _2 =0},$$
however these calculations involve special functions. To avoid this,
we use the problem \eqref{jointLaplaceeq} to write down ODEs for $E[ \tau ^m (x) A^n(x)];$ in the next subsection, we
start with $m=n=1.$
\subsubsection{Correlation between $A(x)$ and $\tau(x)$}
We have:
$$ \frac \partial { \partial \lambda _1} M_\lambda (x) \Big | _ { \lambda_1 = \lambda _2 =0} = -E( \tau (x) ), \
 \frac \partial { \partial \lambda _2} M_\lambda (x) \Big | _ { \lambda_1 = \lambda _2 =0} =  -E( A (x) ), $$
$$ \frac {\partial ^2} { \partial \lambda _1 \partial \lambda _2 } M_\lambda (x) \Big | _ { \lambda_1 = \lambda _2 =0} =
 E( \tau (x) A(x) ),$$
and:
$$ \frac {\partial ^2} { \partial \lambda _1 \partial \lambda _2 } [(\lambda _1 + \lambda _2 x) M_\lambda (x)] \Big | _ { \lambda_1 = \lambda _2 =0} =
-x E( \tau (x)) - E( A(x)). $$
Therefore, by applying $\partial ^2 / \partial \lambda _1 \partial \lambda _2$ to both members of  \eqref{jointLaplaceeq},
calculating for  $\lambda _1 = \lambda _2 =0,$  we obtain:
$$ \frac 1 2 \frac {d^2} {d x^2} E [ \tau (x) A(x)] - \mu \frac d {dx} E [ \tau (x) A(x) ] = -x E( \tau (x)) - E(A(x)).$$
Thus, recalling that $ E( \tau (x)) = x / \mu , \ E(A(x)) = x / 2 \mu ^2 + x^2 / \ 2 \mu ,$ we get:
\begin{Theorem}
$V(x):= E[\tau (x))A(x)] $ is the solution of the problem:
\begin{equation} \label{problemforcorr}
\begin{cases}
\frac 1 2 V''(x) - \mu V' (x) = - \frac { 3 x^2 } { 2 \mu } - \frac x {2 \mu ^2 } , \\ V(0)=0, \\ \lim _ { \mu \rightarrow + \infty } V(x)=0  .
\end{cases}
\end{equation}
\end{Theorem}
\hfill $\Box$

\begin{Remark}
The conditions in $0$ and for $ \mu \rightarrow + \infty $ follow from the facts that,
if the starting point is $0, $ then $\tau (0)=0$ and $A(0)=0, $ while if the drift is $- \infty ,$ the process $X$ reaches
$0$ immediately.
\end{Remark}
By solving  \eqref{problemforcorr}, we obtain:
\begin{equation} \label{EtauA}
V(x)= E[ \tau (x) A(x)] = \frac {x^3} { 2 \mu ^2 } + \frac { x^2 } { \mu ^3 } + \frac x {\mu ^4 } \ .
\end{equation}
The quantity $V(x)$ allows to compute explicitly the correlation
$$ \rho (x) = \frac {cov(\tau (x), A(x)) } {\sqrt { Var( \tau(x)) Var (A(x)) } } =
\frac {E( \tau (x) A(x) ) - E( \tau (x)) E(A(x)) } {\sqrt { Var( \tau(x)) Var (A(x)) } }
 \ , $$
between $\tau(x)$ and $A(x).$ Recalling that $Var ( \tau (x)) = x / \mu ^3, \  Var (A(x)) = x^3 / 3 \mu ^3 + x^2 / \mu ^4 + 5x / \ 4 \mu ^5,$
we finally find:
\begin{equation}
\rho (x) = \left [ \frac {3 x^2 \mu ^2 + 12 \mu x + 12} {4 x^2 \mu ^2 + 12 \mu x +15 } \right ] ^ {1/2} >0 .
\end{equation}
Thus, $A(x)$ and $\tau(x)$ are positively correlated for any $x$ and $\mu .$
By the substitution $\gamma = \mu x,$ one gets:
\begin{equation} \label{exactrho}
\rho (x) =\left ( \frac { 3 \gamma ^2 + 12 \gamma +12} {4 \gamma ^2 +12 \gamma +15 } \right )^ {1/2} ,
\end{equation}
that was derived also in \cite{kearney:jph14} by different arguments.
For $ \mu \rightarrow + \infty $ it holds
$ \rho (x)\rightarrow \sqrt {\frac 3 4 } \simeq 0.866 ,$ while, in the zero drift case $ (\mu \rightarrow 0^+)$ one has
$\rho (x) =  \sqrt {\frac 4 5 } \simeq 0.894 $ (notice that, for $\mu =0 , \ E(A(x))$ and $\tau (x)$ are both infinite).
The knowledge of the exact formula of the correlation $\rho$ allows to verify the consistence of numerical investigation.
Indeed, for several values of $\mu <0, $ we have simulated 100,000 trajectories of the process $X,$ starting e.g. from $x=10,$ and
we have estimated $\tau(x), \ A(x)$ and $\rho (x).$ The results are shown in the Figure 2.

\begin{figure}
\centering
\includegraphics[height=0.3\textheight]{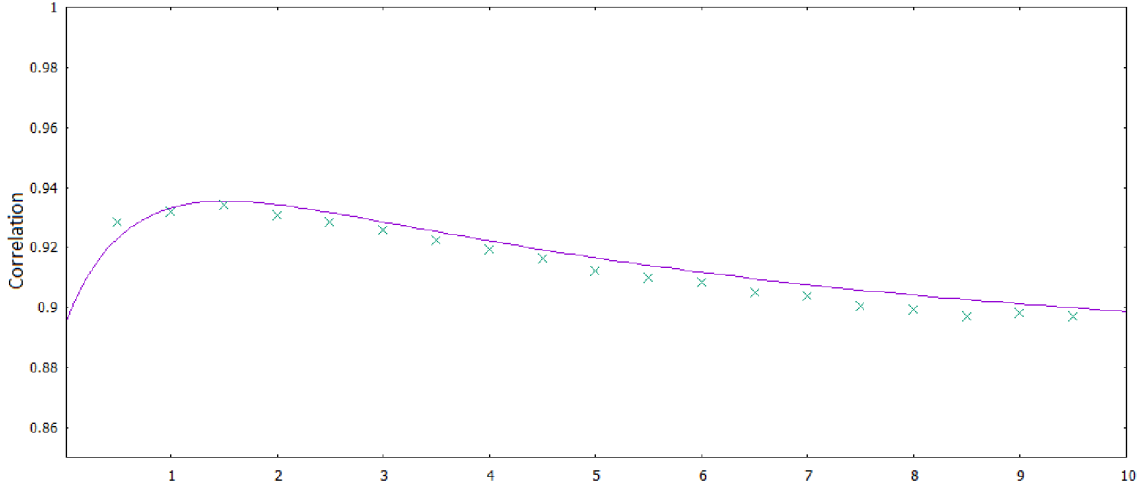}
\caption{Numerical estimate of the correlation $\rho (x)$ between the FPT and the FPA of BM with drift $ - \mu $ starting at $x =10,$ for several values of $\mu. $
On the horizontal axes the values of $\gamma = 10 \mu ;$ the solid line represents the exact curve given by \eqref{exactrho}.
 }
\end{figure}

\subsubsection{Joint moments of any order}
For $m, n \in \mathbb{N},$  we go to study the quantity:
\begin{equation}
V_{m,n}(x): = E \left [\tau(x)^m A(x)^n \right ].
\end{equation}
Of course, $V_{1,1}(x) = E [ \tau(x) A(x) ]$ coincides with the function $V(x)$ already studied.
Applying
$ \frac { \partial ^m } {\partial \lambda _1 ^m} \frac { \partial ^n} {\lambda _2 ^n }$ to both members of  \eqref{jointLaplaceeq},
calculating for  $\lambda _1 = \lambda _2 =0,$
we get that $V_{m,n} (x)$ turns out to be  solution to:
$$ \frac 1 2 \frac {d^2 } {dx^2} \left [ (-1)^ {m+n} V_{m,n} (x) \right ] - \mu \frac d {dx} \left [ (-1)^ {m+n} V_ {m,n} (x) \right ] $$
$$ = m (-1)^{m+n} V _{m-1, n } (x) + nx (-1)^ {m+n} V _ {m, n-1} (x). $$
Thus, rearranging, we  obtain the following:
\begin{Theorem}
$V_{m,n} (x)= E \left [ \tau ^m (x) A^n (x) \right ]$ is the solution of the problem:
\begin{equation} \label{Vmn}
\begin{cases}
\frac 1 2 V''_{m,n} (x) - \mu V' _{m,n} (x) = - m V_{m-1, n} (x) - nx V_{m, n-1} (x) \\ V_{m,n} (0)=0 \\\lim _ { \nu \rightarrow + \infty } V_{m,n} (x) =0
\end{cases} .
\end{equation}
\end{Theorem}

\hfill $\Box$

\noindent Notice that the two conditions are obtained in the same way as in the case of $V_{1,1} (x).$\par\noindent
By solving \eqref{Vmn} e.g. for $m=2, \ n =1,$ one finds:
$$ V_{2, 1} (x)= E \left [ \tau ^2 (x) A(x) \right ] = \frac { x^4} { 2 \mu ^3} + \frac {2 x^3 } {\mu ^4 }
+ \frac {4 x^2 } {\mu ^5 } + \frac { 4x} {\mu ^6 } ,$$
and similar polynomial expressions can be obtained for any $m$ and $n,$ implying that, for $\mu >0,$  the moments
$E \left [ \tau ^m (x) A(x) ^n \right ] $ are finite, for all $m$ amd $n.$ As far as
the form of the solution $V_{m,n}$ of \eqref{Vmn} is concerned, the following holds:

\begin{Proposition} \label{formofVmn}
For all  integers $m,n \ge 0$ the solution of \eqref{Vmn} is a polynomial of degree $m+2n ,$ which vanishes at zero.
\end{Proposition}
{\it Proof.} \ We already know that $V_{1,0} (x) = E( \tau (x)) = x / \mu , \ V_{0,1} (x) = E(A(x)) = x / 2 \mu ^2 + x^2 / 2 \mu ,$ and
$V_{1,1} (x) = E( \tau(x) A(x)) = x^3 / 2 \mu ^2 + x^2 / \mu ^3 + x / \mu ^4 ,$ which are polynomials  vanishing at $x=0,$
with degree, respectively, $1, \ 2$ and $3,$ all corresponding to the value $m+2n.$ As for the quantity $V_{0,n} (x) = E( A(x) ^ n),$
by using \eqref{eqmomentsimple} with $U(x)=x,$ we find that it is the solution of
$$ \begin{cases} \frac 1 2 V'' _{0,n} (x) - \mu V' _{0,n} (x) = - nx V_{0, n-1} (x) \\ V_{0,n} (0)=0 \\ \lim _ { \mu \rightarrow + \infty } V _{0,n } (x) =0 ,
\end{cases}
$$
which is, as easily seen, a polynomial of degree $2n;$ in fact, the homogeneous solution of the ODE is
$c_1 + c_2 e^ { 2 \mu x },$ with $c_i$  constants to be determined, while a particular solution is a polynomial of degree $2n ,$
vanishing at $x=0.$ Since the conditions imply  the constant $c_i$ to be zero, the result follows. Now, let us suppose that the
thesis is true for $m=0$ and $m=1,$ and $n-1$ in place of $n,$  namely we suppose that $V_{0, n-1}$ and $V_{1, n-1}$
are polynomials with degree, respectively, $2(n-1)$ and $2n-1;$  $V_{1,n} (x)$ turns out to be the solution of:
\begin{equation} \label{Vmnbis}
\begin{cases}
\frac 1 2 V''_{1,n} (x) - \mu V' _{1,n} (x) = - m V_{0, n} (x) - nx V_{1, n-1} (x) \\ V_{1,n} (0)=0 \\\lim _ { \nu \rightarrow + \infty } V_{1,n} (x) =0 .
\end{cases}
\end{equation}
The homogeneous solution is the same as the previous ODE, while, being the second member a polynomial of degree
$degr [-mV_{0,n} (x) - x(n-1) V_{1, n-1} (x) ] = 2n,$  the particular solution is a polynomial of degree $2n+1,$ vanishing at zero; hence,
since the imposed conditions imply that the homogeneous part disappears, the solution $V_{1,n} (x)$ is a polynomial of degree
$1+2n,$ which vanishes at zero. Therefore, the thesis is proved for $m=1$ and  every $n.$ \par\noindent
Suppose now the thesis is true for $V_{h,k}(x)$ with $k \in {\bf N}$ if $h < m$ and $k<n$ if $h=m;$ the inductive step consists in proving the thesis for $V_{m,n}.$
But, as in the preceding case, the solution is given by the usual homogeneous solution
(which the conditions require to be zero), plus a polynomial, vanishing at $x=0,$ which has degree $m+2n,$
because the second member of \eqref{Vmn} is a polynomial wit degree
$degr [-m V _ {m-1,n} (x) -  nx V_{m, n-1} (x) ]= m+ 2n -1.$ \par
\hfill $\Box$

\bigskip

\noindent In the remaining part of this subsection we present
a recursive method to find the polynomials $V_{m,n} (x),$ for any $m$ and $n.$ \par\noindent
By using  Proposition \ref{formofVmn}, we obtain that for every $m, \ n \in {\bf N} ,$ there exist real numbers $a_1^{(m,n)}, \dots , a_{m+2n} ^ {(m,n)}$ such that
$$ V_{m,n} (x)= \sum _ {k=1} ^ {m+2n} a_k ^ {(m,n)} x^k ,$$
and, in analogous way, there exist real numbers $a_k^{(m-1,n)}$ for $k=1,  \dots , m+2n -1 $ and $a_k ^ {(m,n-1)}$ for $k=1, \dots , m+2n -2,$
such that
$$ V_{m-1, n} (x)= \sum _ {k=1} ^ {m+2n -1} a_k ^ {(m-1, n)} x^k = a_1 ^ {(m-1, n)} x + \sum _ {k=1} ^ {m+2n -2} a_{k+1} ^ {(m-1,n)} x^ {k+1} ,$$
and
$$V_ {m,n-1} (x)= \sum _ {k=1} ^ {m+2n -2} a_k ^ {(m,n-1)} x^k .$$
By calculating the derivatives involved, the equation in \eqref{Vmn} can be written as
$$ \sum_ {k=1} ^{m+2n-1} k \left [ \frac 1 2 a_{k+1} ^{(m,n)}
(k+1) - \mu a_k ^{(m,n)} \right ] x^{k-1} - \mu a_{m+2n} ^{(m,n)}
(m+2n)x^{m+2n-1} $$
$$= - m a_1 ^{(m-1,n)} x + \sum _{k=1}
^{m+2n-2} (- m a_{k+1} ^{(m-1, n)} - n a_k ^{(m,n-1)} )x^{k+1} .$$
By equaling the coefficients of the two members, we get: \par\noindent
- for constant term,
$$ a_2 ^{(m,n)} - \mu a_1 ^{(m,n)} =0;$$
- for term of degree $1,$
$$ 2 \left ( \frac 3 2 a_3 ^{(m,n)} - \mu a_2 ^{(m,n)} \right ) = - m a_1 ^{(m-1,n)} ;$$
- for term of degree $d=2, \dots, m+2n-2 ,$
$$ (d+1) \left ( \frac {d+2} 2 a_{d+2} ^{(m,n)} - \mu a_{d+1} ^{(m,n)} \right ) = - m a_d ^{(m-1,n)} - n a_{d-1} ^{(m,n-1)} ;$$
- finally, for term of degree $m+2n-1 ,$
$$ - \mu (m+2n) a_{m+2n} ^{(m,n)} = - m a_ {m+2n-1} ^{(m-1,n)} - n a_{m+2n-2} ^{(m,n-1)} .$$
In order to obtain a more compact notation, we introduce the
three matrices \par\noindent $A^{(m,n)}, \ B ^{(m-1,n)}, \ C ^{(m, n-1)} \in \mathbb{R} ^{(m+2n)\times (m+2n)} ,$ such that
$$ A^{(m,n) } _{i,j} =
\begin{cases}
-i \mu   &  i=j \\ \frac 1 2 i (i+1) &   j=i+1, \ i <m+2n \\ 0  & otherwise
\end{cases}
$$
and
$$ B ^{(m-1,n) } = \begin{pmatrix}
0& \dots &  0  \\   . & \dots &  0  \\   0 & \dots  &  -m I_ {m+2n-1}
\end{pmatrix} , \ C ^{(m, n-1)} =
\begin{pmatrix}
0& \dots &  \dots & 0  \\   . & 0 & \dots &  0  \\ . & . & \dots &  0 \\  0 & \dots & .  &  -n I_ {m+2n-1},
\end{pmatrix} ,$$
where $I_k \in \mathbb{R} ^{k\times  k}$ is the identity matrix.
Denoting by $\underline a _{m,n} \in \mathbb{R}^{m+2n}$ the vector of the coefficients of $V _{m,n},$ that is
$$\underline a _{m,n} = \left (a_1 ^{(m,n)}, \dots ,  a_{m+2n} ^{(m,n)} \right ) ^T ,
$$
the following matrix equation resumes the previous equalities:
$$ A^{(m,n)} \underline a _{m,n} = B ^{(m-1,n)} \begin{pmatrix} 0 \\ \underline a _{m-1, n} \end{pmatrix} +
C ^{(m, n-1)}\begin{pmatrix} 0 \\ 0 \\ \underline a _{m, n-1} \end{pmatrix}.$$
Being $A^{(m,n)}$ invertible, it results
\begin{equation} \label{eqforamn}
\underline a _ {m,n} = ( A ^{(m,n)} ) ^{-1} B ^{(m-1,n)} \begin{pmatrix} 0 \\ \underline a _{m-1, n} \end{pmatrix} +
( A ^{(m,n)} ) ^{-1} C ^{(m, n-1)}\begin{pmatrix} 0 \\ 0 \\ \underline a _{m, n-1} \end{pmatrix}.
\end{equation}
As it easily seen, the inverse matrix $( A ^{(m,n)} ) ^{-1}$ is
$$ ( A ^{(m,n)} ) ^{-1} = \begin{pmatrix} -1/ \mu & c_{1,2} & \dots & \dots & \dots & \dots & c_{1, m+2n} \\
0 & -1 / 2 \mu & c_{2,3} & . & \ & \ & . \\
0 & .  & . & . & \ & \ & . \\
0 & . & . & - 1/ i \mu & c_ {i, i+1} & \ & . \\
. & . & . & . & . & . & . \\
. & . & . & . &. & \ & c_{m+2n, m+2n-1} \\
0 & . & . & . &. & 0 & -1 / (m+2n) \mu
\end{pmatrix},
$$
where, for every $i,$
$$c_{i,i+1} = - \frac 1 {2 \mu ^2 },$$
and for every $j > i+1,$
$$c_{i,j} = \frac {(i+1)\cdots (j-1) } {2^{j-i } \mu ^{j-i+1} }.$$
Now, defining the two new matrices
$$ D^{(m-1, n)} := ( A ^{(m,n)} ) ^{-1} B^{(m-1,n)} \ {\rm and } \ E^{(m, n-1)} := ( A ^{(m,n)} ) ^{-1} C^{(m,n-1)},$$
we have
$$ D^{(m-1,n)} _{i,j} = \begin{cases}
0 & j=1 \\ \frac m {i \mu  } & i=j >1 \\ -mc_{i, i+1} = \frac m {2 \mu ^2 } & j=i+1 \\
-m c_{i,j} = m \frac {(i+1) \cdots (j-1) } {2^{j-i} \mu ^{j-i+1} } & j> i+1 \\ 0 & otherwise
\end{cases}
$$
and
$$ E^{(m,n-1)} _{i,j} = \begin{cases}
0 & j=1,2 \\ \frac m {i \mu  } & i=j >2 \\ -nc_{i, i+1} = \frac n {2 \mu ^2 } & j=i+1, j>2 \\
-n c_{i,j} = n \frac {(i+1) \cdots (j-1) } {2^{j-i} \mu ^{j-i+1} } & j> i+1, j >2 \\ 0 & otherwise
\end{cases}
$$
Finally, \eqref{eqforamn} becomes
\begin{equation} \label{iterative}
\underline a _ {m,n} = D ^{(m-1, n)} \begin{pmatrix} 0 \\ \underline a _{m-1, n} \end{pmatrix} +
E ^{(m, n-1)} \begin{pmatrix} 0 \\ 0 \\ \underline a _{m, n-1} \end{pmatrix},
\end{equation}
which provides a recursive formula to find the coefficients $a_k^{(m,n)}, \ k=1, \dots, m+2n.$ Thanks to the fact that the involved matrices are
triangular, for $m$ and $n$ not too large, \eqref{iterative} represents a faster way to obtain  the coefficients of the polynomial
$V_{m,n} (x),$ than solving directly
the ODE \eqref{Vmn}.
\subsection{Expected value of the time average till the FPT}
In this subsection, we fill find a closed form for
$$E \left ( \frac {A(x)} { \tau (x)}  \right )= E \left ( \frac 1  { \tau (x)} \int _0 ^{ \tau (x)} X(t) dt  \right ),$$
that is the expected value of the time average of drifted BM $X(t)= x - \mu t + B_t $ till its FPT below zero. \par\noindent
If $M_ \lambda (x)= M_ {\lambda _1, \lambda _2} (x)$ is the joint Laplace transform of $(A(x), \tau (x)) ,$
defined by  \eqref{joint Laplace},
we notice that
$$E \left ( A(x) e^{\lambda _1 \tau (x)} \right ) = - \lim _{ \lambda _2 \rightarrow 0^+ } \frac \partial { \partial \lambda _2 }
M_ {\lambda _1, \lambda _2} (x),$$
and
\begin{equation} \label{A/tau}
E \left ( \frac {A(x)} { \tau (x)}  \right ) = \int _0 ^{ + \infty } E \left ( A(x) e^{\lambda _1 \tau (x)} \right ) d \lambda _1 \ .
\end{equation}
Thus, the calculation of $E(A(x)/  \tau(x) )$ is reduced to find
$- \lim _{ \lambda _2 \rightarrow 0^+ } \frac \partial { \partial \lambda _2 }
M_ {\lambda _1, \lambda _2} (x).$
By taking the partial derivative with respect to $\lambda _2$ in \eqref{jointLaplaceeq},  calculated for $\lambda _2 =0,$ we easily obtain that
$w(x):= E( A(x) e^{- \lambda _1 \tau (x)} )$ is solution to the equation
$$ \frac 1 2 w'' (x) - \mu w'(x)=   \lambda _1 w(x)  -x E ( e^ {- \lambda _1 \tau (x)}).$$
By recalling the expression of the Laplace transform of $\tau (x),$ that is \par\noindent
$E ( e^ {- \lambda _1 \tau (x)})= e^{ - x ( \sqrt {\mu ^2 + 2 \lambda _1 } - \mu ) }$ (see  \eqref{Laplacetau}), and reasoning as before,
we obtain that $w(x)$ is the solution
of the problem
\begin{equation}
\begin{cases}
\frac 1 2 w'' (x) - \mu w'(x) =  \lambda _1 w(x)  - x e^{ - x ( \sqrt {\mu ^2 + 2 \lambda _1 } - \mu ) } \\
w(0)=0 \\ \lim _{ \mu \rightarrow + \infty } w(x) =0 .
\end{cases}
\end{equation}
By calculation, we find that the explicit expression of $w(x)$ is:
\begin{equation}
w(x)= E \left [A(x) e^{- \lambda _1 \tau (x)} \right ] = e^{ \mu x } e^{- x \sqrt {\mu ^2 + 2 \lambda _1 } }\left [ \frac x {2 (\mu ^2 + 2 \lambda _1 ) } +
\frac {x^2 } { 2 \sqrt { \mu ^2 + 2 \lambda _1 } } \right ] ,
\end{equation}
which coincides with equation (12) of \cite{kearney:jph14}, that was derived by using different arguments.
Finally, integrating with respect to $\lambda _1,$ from \eqref{A/tau} we get that the expected value of the time average of $X$ till its FPT below zero is:
\begin{equation}
E[A(x) / \tau  (x)] = \frac x 2 \left [ 1 + \int _0 ^{ + \infty } \frac {e^{-sx} } {s + \mu } \ ds  \right ] .
\end{equation}
This result was already obtained in \cite{kearney:jph14} (see eq. (29) therein), as a consequence of \eqref{A/tau}.
Notice that $E[A(x) / \tau  (x)]$ turns out to be $\ge x/2,$ being $\lim _{ \mu \rightarrow + \infty} E[A(x) / \tau  (x)] =x/2,$ and
$\lim _ { \mu \rightarrow 0^+} E[A(x) /  \tau  (x)] = + \infty .$

\section{Conclusions and final remarks}
In this paper,
we have continued the study,  already undertook in \cite{abundo:mcap13},
of the first-passage area (FPA)  $A(x),$ swept out by a one-dimensional jump-diffusion process $X(t),$ starting from $x>0,$
till its first-passage time (FPT) $\tau(x)$
below zero.
Here, we have investigated the joint distribution of $\tau(x)$ and $A(x),$ in the special case when
$X$ is Brownian motion with negative drift $- \mu,$ that is, $X(t)= x - \mu t + B_t, $ ($B_t$ denotes standard Brownian
motion).
We have established differential equations with boundary conditions for the Laplace transform
of the random vector $(\tau(x), A(x)),$ and for
the joint moments  $E[\tau(x)^m A(x)^n]$ of the FPT and FPA; moreover, we have presented an  algorithm  to find recursively
$E[\tau(x)^m A(x)^n],$ for any $m$ and $n.$
In particular, closed formulae for $E[\tau(x) A(x)]$ and then for the correlation between $\tau(x)$ and $A(x),$ were
explicitly obtained, which match existing results, obtained by approximation arguments in \cite{kearney:jph14}.
Furthermore, we have calculated the expected value of the time average of $X(t)$ till its FPT below zero. \par\noindent
The difference compared to similar articles (see e.g. \cite{kearney:jph14}) is that our results have been obtained
without the use of special functions. \par
We remark that all the results contained in this paper can be extended in principle to the case of a one-dimensional time-homogeneous
diffusion $X(t),$
which is the solution of a stochastic differential of the form
\begin{equation} \label{diffueq}
\begin{cases}
dX(t) = b(X(t))dt +\sigma(X(t)) dB_t \\
X(0)=x >0 ,
\end{cases}
\end{equation}
where the functions $b(\cdot), \sigma (\cdot)$
satisfy suitable conditions for the
existence and  uniqueness of the solution of \eqref{diffueq} (see e.g. \cite{gimsko:sde72}, \cite{ikwa:sde81}).
Now, the FPT of $X$ below zero is defined by $\tau (x)= \inf \{t>0: X(t) \le 0  | X(0)=x \}$ and it is required the additional condition
that $\tau(x)$ is finite with probability one (this evidently depends on the behavior of the drift $b(x)$ and the
diffusion coefficient $\sigma (x)).$
Under this condition, all the differential problems concerning the various quantities
considered in Section 2 still hold, provided that  $\frac 1 2 V''(x) - \mu V'(x)$ is replaced with the infinitesimal operator
of $X,$ i.e. $LV(x)= \frac 1 2 \sigma ^2 (x) V''(x) + b(x)V'(x).$ Of course, the differential equations obtained in this way are
more complicated, since in  most cases their solutions cannot be obtained in closed form.


\begin{thebibliography}{spc}

\bibitem [1] {abundo:mcap13}
Abundo, M. \newblock
On the first-passage area of a one-dimensional jump-diffusion process.
\newblock{Methodol Comput Appl Probab } 2013, 15, 85-–103.
(Online First, April 2011). DOI 10.1007/s11009-011-9223-1.


\bibitem[2]
{gimsko:sde72}
Gihman, I.I. and  Skorohod, A.V.\newblock
Stochastic differential
equations. \newblock
Springer-Verlag, Berlin, 1972.

\bibitem[3]
{grand:tab80}
Grandshteyn, I.S. and Ryzhik, I. M. \newblock
Tables of Integrals, Series and Products. \newblock
5th ed. Academic, London, 1980.



\bibitem[4]
{ikwa:sde81}
Ikeda, N. and Watanabe, S.\newblock
Stochastic differential equations and
diffusion processes. \newblock
North-Holland Publishing Company, 1981

\bibitem[5]
{janson:pro07}
Janson, S.\newblock
Brownian excursion area, Wright's
constants in graph enumeration, and
other Brownian areas. \newblock {Probability Surveys} 2007,
     4: 80--145.

\bibitem[6]
{kartay:sto75}
Karlin, S. and Taylor, H.M.\newblock
A second course in stochastic processes. \newblock
Academic Press, New York, 1975.

\bibitem[7]
{kearney:jph14}
Kearney, M.J., Pye, A.J., and Martin R.J. \newblock
On correlations between certain random variables associated with first passage
Brownian motion. \newblock
{J. Phys. A: Math and Theor.} 2014, 47 (22): 225002. DOI: 10.1088/1751-8113/47/22/225002

\bibitem[8]
{kearney:jph05}
Kearney, M. J. and Majumdar, S.N.\newblock
On the area under a continuous time Brownian motion till its first-passage time. \newblock {J. Phys. A: Math. Gen.} 2005,
     38: 4097--4104.

\bibitem[9]
{kearney:jph07}
Kearney, M. J., Majumdar, S.N. and Martin R.J.\newblock
The first-passage area for drifted Brownian motion and the moments of the Airy distribution.
\newblock {J. Phys. A: Math. Theor.} 2007,
     40: F863--F864.




\bibitem[10]
{knight:jam00}
Knight, F. B. \newblock
The moments of the area under reflected Brownian Bridge conditional on its local time at zero.
\newblock {Journal of Applied Mathematics and Stochastic Analysis} 2000,
     13 (2): 99--124.




\bibitem[11]
{perman:aap96}
Perman, M. and Wellner, J. A. \newblock
On the Distribution of Brownian Areas. \newblock
{The Annals of Applied Probability} 1996, 6 (4) 1091--1111.




\end{thebibliography}
\end{document}